\documentclass[brochure,10pt]{article}
\usepackage{amssymb,amsfonts,amsmath,footnote}
\usepackage{latexsym}
\usepackage{amsmath}
\usepackage{amsthm}
\usepackage{amsfonts}
\usepackage{amssymb}
\usepackage{amsmath,amscd}
\usepackage{graphicx}
\usepackage{graphics}
\usepackage{etex}
\usepackage[matrix,arrow]{xy}
\usepackage{color}
\usepackage{tikz-cd}
\newtheorem{Th}{Theorem}

\newtheorem{Lm}{Lemma}

\newcommand{\be}{\begin{equation}}
\newcommand{\ee}{\end{equation}}
\newcommand{\bes}{\begin{equation*}}
\newcommand{\ees}{\end{equation*}}

\newcommand{\R}{\mathbb{R}}

\newcommand{\C}{\mathbb{C}}
\newcommand{\Z}{\mathbb{Z}}

\def\La{\Lambda}
\def\La{\Lambda}

\def\lf{\left}
\def\rg{\right}

\def\al{\alpha}
\def\la{\lambda}

\def\ep{\varepsilon}

\def\ds{\displaystyle}

\def\om{\omega}
\def\p{\partial}

\def\r{\mathfrak r}

\begin{document}

\title{Almost Monotonicity Formula for H-minimal Legendrian Surfaces in the Heisenberg Group. }

\author{ Tristan Rivi\`ere\footnote{Department of Mathematics, ETH ,
Switzerland.}}

\maketitle

{\bf Abstract :}{\it We prove an almost monotonicity formula for H-minimal  Legendrian Surfaces (also called {\it contact stationary legendrian immersions}) in the Heisenberg Group ${\mathbb H}^2$ . From this formula we deduce a Bernstein-Liouville type theorem
for   H-minimal  Legendrian Surfaces. This last result happens to be in particular a crucial ingredient in the analysis part of the proof of the Willmore conjecture given by the author in \cite{Riv}. }

\medskip

\noindent{\bf Math. Class. 53D10, 53D12, 53A10, 49Q05 }

\medskip

\noindent{\bf  Key words.} {\it Hamiltonian-minimal legendrian Surfaces, contact stationary legendrian surfaces,  stationary exact lagrangian surfaces,  monotonicity formula, Heisenberg group, Bernstein theorem, Liouville theorem.}

\section{Introduction}
Conserved quantities and monotonicity formula are fundamental notions in the calculus of variations.  These identities are reflecting the existence of groups of symmetries of the underlying lagrangian most of the time in locally isotropic spaces. The most illustrative exemple is  maybe the monotonicity formula for minimal immersion of a surfaces in an euclidian space ${\R}^n$. This identity says the  following: let  $\Phi$ be a minimal immersion of a surface $\Sigma$ without boundary into ${\R}^n$, then the following {\bf conservation law} holds
\be
\label{001}
\forall\, r>0\quad\quad\frac{1}{r^2}\int_{\rho<r}\ dvol_\Sigma=\int_{\rho<r}\frac{|(\nabla\rho)^\perp|^2}{\rho^2}\ dvol_\Sigma+\pi\,  \mbox{Card}\lf( \Phi^{-1}(0)\rg)\ ,
\ee
where $\rho$ is the distance function to the origin of the space,  $\ dvol_\Sigma$ denotes the volume form on $\Sigma$ induced  by the immersion $\Phi$ , $(\nabla^\Sigma\rho)^\perp$ is the projection to the normal codimension 2 plane to the immersion of the gradient  of $\rho$ and $\mbox{Card}\lf( \Phi^{-1}(0)\rg)$ is the number of pre-images of the origin 0 by the immersion $\Phi$ (see for instance \cite{Sim}). This  formula which holds as well for  much weaker mathematical objects critical point of the area than minimal immersions (such as {\it stationary varifolds} for instance) is  the starting point for the analysis of the variation of the area in euclidian or even riemannian spaces.

The purpose of the present work is to look for the natural replacement of (\ref{001}) in a ``fully anisotropic environnement'' that is in an ambient space where all directions are not equal (even infinitesimally). The most symmetric (or isotropic) model for  anisotropy is maybe given by the Heisenberg groups ${\mathbb H}^n$ of dimension $2n+1$ in which one directions is ``forbidden'' while perfect isotropy holds in the $2n$ remaining ones. 
The coordinates in ${\mathbb H}^n$ will be denoted $(z_1,\cdots z_{2n},\varphi)$ where the last coordinate $\varphi$ is called the {\it legendrian coordinate}. The so called {\it Horizontal Hyperplanes} $H$ are generated at every points by the following $2n$ vectors
\[
X_i:=\frac{\p}{\p z_{2i-1}}-z_{2i}\frac{\p}{\p\varphi}\ \quad,\quad\ Y_i:=\frac{\p}{\p z_{2i}}+z_{2i-1}\frac{\p}{\p\varphi}\quad\mbox{ for }i=1\cdots n
\]
We take these vectors to realize an orthonormal frame in such a way that the canonical projection $\pi$ from $T{\mathbb H}^{n}$ into $T{\mathbb R}^{2n}$ given by
\[
\pi_\ast X_i=\frac{\p}{\p z_{2i-1}}\quad,\quad \pi_\ast Y_i=\frac{\p}{\p z_{2i}}\quad\mbox{ and }\quad \pi_\ast\frac{\p}{\p\varphi}=0
\]
realizes at every point an isometry from $H$ into $T{\mathbb R}^{2n}$. An immersion $\La$ of a surface $\Sigma$ into ${\mathbb H}^n$ is called {\it legendrian}  if it is tangent to $H$ at every point. This is also equivalent to the following {\it contact condition}
\[
\La^\ast\al=0\quad\mbox{ where}\quad\al:=-d\varphi +\sum_{i=1}^n z_{2i-1}\,dz_{2i}-z_{2i}\,dz_{2i-1}\quad.
\]
The projection by $\pi$ of a {\it legendrian immersion}  $G:=\pi\circ\La$ generates obviously\footnote{ The fact that $d\pi$ is an isometry from $H$ into $T{\mathbb R}^2n$ and that $\La$ is an immersion everywhere tangent to $H$ preserves the fact that $G:=\pi\circ\La$ is an immersion, moreover $d\al=\pi^\ast\om$ and $\La^\ast\al=0$ gives $\La^\ast\pi^\ast\om=G^\ast\om=0$. }a {\it Lagrangian immersion} of ${\R}^{2n}$ that is an immersion into ${\R}^{2n}$ satisfying
\[
G^\ast\om=0\quad\mbox{ where}\quad \om:=2\,\sum_{i=1}^n dz_{2i-1}\wedge dz_{2i}\ .
\]
The reciproque holds as well locally but we prefer to work with Legendrian immersion rather than Lagrangian immersion because the \underbar{global} existence of the legendrian coordinates is necessary for  monotonicity formula to hold for the area functional  as observed first by R.Schoen and J.Wolfson in \cite{SW1} (see  a counter-exemple pages 192-193). Lagrangian surfaces with a global legendrian lift are also called {\it exact lagrangian surfaces}.

Critical points to the area functional for arbitrary compactly supported perturbations among legendrian immersions are called {\it hamiltonian stationary legendrian immersions} or simply {\it  H-minimal Legendrian Surfaces}. Since $d\pi$ is an isometry from $H$ into $T{\mathbb R}^{2n}$  and since locally a one to one correspondence holds between legendrian and Lagrangian immersions, such an immersion is projected by $\pi$ onto a critical point of the area among lagrangian surfaces. These surfaces have been considered first by Oh under the name of 
{\it Hamiltonian Stationary Surfaces}  (see \cite{Oh1} and \cite{Oh2}).  For such a surface the mean-curvature vector is given by the image by the canonical complex  structure $i$ of ${\R}^{2n}\simeq {\C}^n$ of the gradient  along the surface of an  {\it harmonic $2\pi{\Z}$ multivalued function } $\beta$, the {\it lagrangian angle}, whose differential is generating the {\it Maslov class} of the lagrangian surface (see also \cite{SW1}). We shall now concentrate  mostly on the case $n=2$ . In this particular dimension Lagrangian surfaces are characterized by the fact that the action of $i$ realizes  an isometry between the tangent 2-space and the orthogonal 2-space to the surface. In conformal parametrization a lagrangian immersion is stationary  if and only if the multivalued function $\beta$ (well defined up to a multiple of $2\pi$) satisfies\footnote{ The Euler Lagrange equation (\ref{equation}) can also be formulated as follows : there exists an $S^1-$valued harmonic map $u$ on $\Sigma$ equipped with the induced metric by the immersion $G$ such that
\[
div^\Sigma(u\,\nabla^\Sigma G)=0\quad.
\]}
\be
\label{equation}
\lf\{
\begin{array}{l}
\ds i\,\Delta G=\, 2\ \nabla\beta\cdot\nabla G\\[5mm]
\ds \Delta\beta=0
\end{array}
\rg.
\ee
An important quantity in the Heisenberg group is given by $\r:=(\rho^4+4\,\varphi^2)^{1/4}$ where $\rho$ denotes the euclidian distance in ${\R}^{4}$ to the origin. The function $\r$ is called the {\it Folland-Kor\'anyi gauge}. It defines the left invariant homogeneous {\it Heisenberg distance}  equivalent to the {\it Carnot Carath\'eodory distance} related to
the minimal length among horizontal geodesics between two points. The striking fact for the interpretation of our main result in relation to more classical monotonicity formula is that in ${\mathbb H}^2$ the Green function for the horizontal laplacian $\Delta_H$ is proportional to $\r^{-2}$.

\medskip

Our main result in the present work is the following result.
\begin{Th}
\label{th-main} {\bf [Almost Monotonicity]}
There exists a universal constant $C>0$ such that for any  smooth H-minimal legendrian immersion $\La$ of an oriented surface $\Sigma$ into ${\mathbb H}^2$.  
Then we have
\be
\label{000}
\forall\, r<1\quad\quad C^{-1}\lf[\theta_0+\int_{\r<r/2} \frac{|(\nabla^\Sigma\r)^\perp|^2}{\r^2} \ dvol_\Sigma\rg]\le \frac{1}{r^2}\int_{\r<r}\ dvol_\Sigma\le C\ \int_{1/2<\r<2}\ dvol_\Sigma\ ,
\ee
where $\theta_0:=2\pi\,\mbox{Card}\,{{\La}}^{\,-1}(0)$ where $\r:=(\rho^4+4\,\varphi^2)^{1/4}$ is the Folland-Kor\'anyi gauge and where $(\nabla^\Sigma\r)^\perp$ denotes the projection of the gradient of the Folland-Kor\'anyi gauge onto the orthogonal 2-plane to the tangent space of the immersion $\La$ within the horizontal plane $H$.
\hfill $\Box$
\end{Th}

\medskip

The Inequality (\ref{000}) is the counterpart for H-minimal legendrian immersion of the monotonicity identity (\ref{001}) for minimal surfaces in ${\R}^n$. In fact inequality (\ref{000}) is the consequence of the {\bf conservation law} (\ref{k-10-rep})  for H-minimal legendrian immersion which corresponds to the conservation law (\ref{k-26-rep}) for classical minimal surfaces.

The control of the area density from above and below given by (\ref{th-main}) is the starting point for the blow-up analysis of $H-$minimal legendrian surfaces in ${\mathbb H}^2$ and is centrally used in the implementation of the viscosity method for  Stationary Lagrangian Surfaces in K\"ahler-Einstein Manifolds in \cite{Piga-Riv-Lag-Visc} which is itself needed for the proof of the {\it Willmore conjecture} given by the author in \cite{Riv}. We expect  the result also to hold in higher dimension and  for much weaker objects than H-minimal immersions such as the weakly conformal exact lagrangian stationary maps in \cite{SW1} and \cite{SW2} or even to the notion of {\it integer rectifiable legendrian contact stationary varifold} mentioned in \cite{SW1} . A control from above and bellow of the area density for exact stationary lagrangian maps of surfaces  has first been given by Schoen and Wolfson in \cite{SW1}, (proposition 3.2)
\footnote{Observe that in the proof of Proposition 3.2 given in \cite{SW2} the computation of $div_0(X_{\eta(t,\theta)})$ at the bottom of page 14 holds away from the origin. For a smooth $H-$minimal legendrian immersion for instance the asymptotic expansion of the  vector field $X_{\theta}$ generated by the hamiltonian $\theta$ and whose flow is preserving the contact distribution is of the form
\[
X_{\theta}=-\frac{1}{\rho}\frac{\partial }{\partial \rho}+O(1)
\]  
where the quantity $\rho^2:=x^2+y^2$ corresponds to $\sqrt{s}$ in \cite{SW2}. Hence the divergence of $X_\theta$ generates a Dirac mass at the each of the pre-images of the origin.  
The  weights in front of these Dirac masses have all the same sign (their value is in fact $-2\,\pi$ in the smooth setting but might vary if the H-minimal legendrian is not exactly an immersion and has {\it Schoen-Wolfson $p-q$ type conical singularities} then the weights should be  $-2\,\pi \sqrt{p\,q}$ .) and are  ``responsible'' for the lower bound  in the area density control.}. 

A consequence of the {\bf conservation law } (\ref{k-10-rep}) is the following Bernstein-Liouville type theorem for  H-minimal legendrian immersion which says roughly that if such an immersion is asymptoticaly
a {\it Lagrangian plane} then it must be a lagrangian plane.
\begin{Th}
\label{th-berns-liouv}{\bf [Bernstein-Liouville type theorem]}
Let $\Lambda$ be a  smooth H-minimal legendrian immersion of an oriented surface $\Sigma$ into ${\mathbb H}^2$. Assume $\La$ is asymptotically a plane at infinity that is
\be
\label{plane}
\lim_{r\rightarrow+\infty}\frac{1}{r}\int_{r<\r<2r} \frac{1}{\r}\ dvol_\Sigma=2\pi\quad\mbox{ and }\quad \frac{\rho}{\r}\rightarrow 1\quad\mbox{ as }\ \r\rightarrow +\infty\ ,
\ee
then $\La(\Sigma)$ is a lagrangian plane in ${\C}^2\simeq {\mathbb H}^2\cap\{\sigma=0 \}$.
\hfill $\Box$
\end{Th}

\section{Preliminaries}
Recall that the Lagrangian projection ${G}:=\pi\circ\La$ from ${\mathbb H}^2$ into ${\C}^2$ of a conformal Hamiltonian Stationary Legendrian Immersion $\vec{\La}$ satisfies locally in conformal charts (see \cite{SW2})
\[
i\,\Delta G=2\, \nabla \beta\cdot\nabla G
\]
Taking the scalar product respectively with $G$ and $iG$ gives
\be
\label{heat-12}
\lf\{
\begin{array}{l}
\ds div(G\cdot i\nabla G)=\nabla\beta\cdot\nabla |G|^2\\[5mm]
\ds div(G\cdot \nabla G)-|\nabla G|^2=-2\, \nabla \beta\cdot (G\cdot i\nabla G)
\end{array}
\rg.
\ee
Since $G$ is lagrangian we have that the multiplication by $i$ realizes an isometry between the tangent and the normal planes to the immersion $G$. Hence
\[
div(G\cdot i\nabla^\perp G)=0\ .
\]
Hence we introduce the legendrian coordinate equal to zero at the origin of ${\mathbb H}^2$ and given by
\[
d\varphi:=G\cdot i dG
\]
We denote $\rho^2:=|G|^2$. With this notation (\ref{heat-12}) becomes (independently of coordinates)
\be
\label{kora-4re}
\lf\{
\begin{array}{l}
\ds \nabla^\Sigma\beta\cdot\nabla^\Sigma\rho^2=\Delta^\Sigma\varphi\\[5mm]
\ds \nabla^\Sigma\beta\cdot\nabla^\Sigma\varphi=1-4^{-1}\,\Delta^\Sigma\rho^2\ .
\end{array}
\rg.
\ee
where $\nabla^\Sigma$ and $\Delta^\Sigma$ denote respectively the {\it gradient operator} and  the {\it negative laplace Beltrami operator} for the induced metric by $G$ (or $\La$) on $\Sigma$. 

\medskip

The {\it lagrangian angle multivalued function} $\beta$ plays the role of a lagrangian multiplier and satisfies the Euler Lagrange equation
\be
\label{heat-12-a}
\Delta^\Sigma\beta=0\ .
\ee

\medskip

Since $(e^{-\la}\p_{x_1}G, e^{-\la}\p_{x_2}G, i\,e^{-\la}\p_{x_1}G, i\,e^{-\la}\p_{x_2}G)$ realizes an orthonormal frame of $G^\ast T{\C}^2$ we have
\[
\rho^2=|G|^2=e^{-2\la}\,|G\cdot \nabla G|^2+e^{-2\la}\,|G\cdot i\,\nabla G|^2=|2^{-1} \nabla^\Sigma\rho^2|^2+|\nabla^\Sigma\varphi|^2
\]
from which we deduce
\be
\label{heat-13}
1=|\nabla^\Sigma\rho|^2+\rho^{-2}\,|\nabla^\Sigma\varphi|^2\ .
\ee

\section{Density Bound and Dirichlet Energy bound of the Arctangent of the Phase. }
We call the phase the scaling invariant quantity in the Heisenberg group given by
\[
\sigma:=\frac{2\varphi}{\rho^2}\ .
\]
The purpose of the present section is to prove that, for any  H-minimal legendrian  immersion,  the sum of the density at the origin with the Dirichlet energy of $\arctan\sigma$ in a unit ball for the Folland-Kor\'anyi distance is controlled by a constant times the area
contained in the surrounding dyadic annulus. Precisely we have
\begin{Lm}
\label{lm-mono}
There exists a universal constant $C_0>0$ such that for any ${\La}$  smooth $H-$minimal legendrian  immersion of an oriented surface $\Sigma$ into ${\mathbb H}^2$. Denoting $\sigma:=2\,\varphi/\rho^2$ 
then we have
\be
\label{k-15-rep}
2\pi\,\mbox{Card}\,{{\La}}^{\,-1}(0)+\int_{\r<1}\lf|\frac{\nabla^\Sigma\sigma}{1+\sigma^2}\rg|^2\ dvol_\Sigma\le C_0\ \int_{1<\r<2}\ dvol_\Sigma\ .
\ee
\hfill $\Box$
\end{Lm}
\noindent{\bf Proof of lemma~\ref{lm-mono}.}
Using (\ref{kora-4re}) away from $\rho= 0$
\be
\label{k-6}
\begin{array}{rl}
\ds\nabla^\Sigma\sigma\cdot\nabla^\Sigma\beta&\ds=\lf[2\, \rho^{-2}\nabla^\Sigma\varphi-2\,\varphi\,\rho^{-4}\,\nabla^\Sigma\rho^2\rg]\cdot\nabla^\Sigma\beta\\[5mm]
 &\ds=\frac{1}{2\,\rho^2}\, \lf[ 4-\Delta^\Sigma\rho^2  \rg] -\frac{2\,\varphi}{\rho^4}\ \Delta^\Sigma\varphi\\[5mm]
 &\ds=2\, \rho^{-2}-2^{-1}\rho^{-2}\,\Delta^\Sigma\rho^2-\rho^{-4} \, \Delta^\Sigma\varphi^2+2\,\rho^{-4}\ |\nabla^\Sigma\varphi|^2\\[5mm]
  &\ds=2\, \rho^{-2}-4^{-1}\rho^{-4}\,2\, \rho^2\,\Delta^\Sigma\rho^2-\rho^{-4} \, \Delta^\Sigma\varphi^2+2\,\rho^{-4}\ |\nabla^\Sigma\varphi|^2\\[5mm]
  &\ds=2\, \rho^{-2}-4^{-1}\rho^{-4}\, \lf[\Delta^\Sigma\rho^4-2\,|\nabla^\Sigma\rho^2|^2+4\,\Delta^\Sigma\varphi^2\rg]+2\,\rho^{-4}\ |\nabla^\Sigma\varphi|^2\\[5mm]
 &\ds=2\, \rho^{-2}+2\,\rho^{-2}\,|\nabla^\Sigma\rho|^2+2\,\rho^{-4}\ |\nabla^\Sigma\varphi|^2-4^{-1}\,\rho^{-4}\,\Delta^\Sigma\r^4\\[5mm]
 &\ds=4\,\rho^{-2}-4^{-1}\,\rho^{-4}\,\Delta^\Sigma\r^4=4\,\frac{\sqrt{1+\sigma^2}}{\r^2}-4^{-1}\,\frac{1+\sigma^2}{\r^4}\,\Delta^\Sigma\r^4
\end{array}
\ee
where we used (\ref{heat-13}). This implies assuming $\rho\ne0$ and dividing by $1+\sigma^2$
\be
\label{k-6-rep}
\begin{array}{rcl}
\ds\frac{\nabla^\Sigma\sigma}{1+\sigma^2}\cdot\nabla^\Sigma\beta&=&\ds\frac{4}{\r^2}\,\frac{1}{\sqrt{1+\sigma^2}}-\frac{1}{4\,\r^4}\,\Delta^\Sigma\r^4\\[5mm]
\ds&=&\ds\frac{4}{\r^2}\,\frac{1}{\sqrt{1+\sigma^2}}-\frac{1}{4}\,\mbox{div}^\Sigma\lf[\frac{\nabla^\Sigma\r^4}{\r^4}\rg]-4\, |\nabla\log\r|^2\\[5mm]
\ds&=&\ds\frac{4}{\r^2}\,\frac{1}{\sqrt{1+\sigma^2}}-4\, |\nabla^\Sigma\log\r|^2-\Delta^\Sigma\log\r
\end{array}
\ee
Since 
\[
\lf\{
\begin{array}{l}
\ds d\varphi=z_2\,dz_1-z_1\,dz_2+z_4\,dz_3-z_3\,dz_4\\[5mm]
\ds\rho\, d\rho=z_1\,dz_1+z_2\,dz_2+z_3\,dz_3+z_4\,dz_4
\end{array}
\rg.
\]
This implies
\be
\label{kora-12a}
|\nabla^H\varphi|=\rho\quad,\quad|\nabla^H\rho|=1\quad\mbox{ and }\quad\nabla^H\varphi\cdot\nabla^H\rho=0\quad.
\ee
We deduce the length of the {\it horizontal gradient} of the {\it Folland - Kor\'anyi gauge}
\be
\label{kora-13}
\ds|\nabla^H\r|^2=   \frac{\rho^6}{\r^6}\, |\nabla^H\rho|^2+\frac{|\nabla^H\varphi^2|^2}{\r^6}= \frac{\rho^2}{\r^2}\, \frac{\rho^4+4\,\varphi^2}{\r^4}= \frac{\rho^2}{\r^2}\ .
\ee
hence we have
\be
\label{kora-13-a}
|\nabla^H\r|^2=\frac{\rho^2}{\r^2}=\frac{1}{\sqrt{1+\sigma^2}}\ .
\ee
Combining (\ref{k-6-rep}) and (\ref{kora-13-a}) gives away from $\rho=0$, since $\Delta^\Sigma\beta=0$,
\be
\label{k-6-rep}
\begin{array}{rcl}
\ds\mbox{div}^\Sigma\lf( \arctan\sigma\ \nabla^\Sigma\beta+\nabla^\Sigma\log\r\rg)&=&\ds{4}\,\frac{|(\nabla^\Sigma\r)^\perp|^2}{\r^2}\ .
\end{array}
\ee
Assume now $\rho=0$ only happens at the points $p$ where $\r(p)=0$ , then, since we have a smooth legendrian immersion we deduce that
\be
\label{k-6-rep}
\begin{array}{rcl}
\ds\mbox{div}^\Sigma\lf( \arctan\sigma\ \nabla^\Sigma\beta+\nabla^\Sigma\log\r\rg)&=&\ds{4}\,\frac{|(\nabla^\Sigma\r)^\perp|^2}{\r^2}+\sum_{\r(p)=0}\theta_0(p)\ \delta_{p}\ .
\end{array}
\ee
where, for $\ep>0$ small enough and fixed
\[
\theta_0(p):=\lim_{t\rightarrow 0}\int_{\{\r=t\}\cap B_\ep(p) \ }\frac{\p_\nu\r}{\r}\ dl_\Sigma\ .
\]
where $\nu$ is the unit normal tangent to $\Sigma$ pointing out of the set  $\r\le t$. In conformal coordinates $x=(x_1,x_2)$ for the associated lagrangian immersion $G$ assuming that $\r\circ{\Lambda}(0)=0$ we have the existence of an orthonormal basis of ${\C}^2$ such that
$J(\vec{\ep}_1)=\vec{\ep}_2$ and $J(\vec{\ep}_3)=\vec{\ep}_4$ and
\[
G(x)=e^{\la_0}\ (x_1\,\vec{\ep}_1+x_2\,\vec{\ep}_3)+O(|x|^2)
\]
This gives 
\[
{G}^{\,\ast}d\varphi={G}^{\,\ast}(z_2\,dz_1-z_1\, dz_2+z_4\, dz_3-z_3\,dz_4)=O(|x|^2)
\]
Hence
\be
\label{k-6-rep-bis}
\varphi(x)=O(|x|^3)\quad, \quad \frac{\rho}{\r}=1+o(1)\quad\mbox{ and }\quad\nabla^\Sigma\r=\nabla^\Sigma\rho\,(1+o(1))
\ee
This gives
\be
\label{k-6-rep-ter}
\theta_0(p):=\lim_{t\rightarrow 0}\int_{\{\r=t\}\cap B_\ep(p) }\frac{\p_\nu\r}{\r}\ dl_\Sigma=\lim_{t\rightarrow 0}\int_{\rho=t}\frac{\p_\nu\rho}{\rho}\ dl_\Sigma=2\pi\,
\ee
We claim that the following holds on the whole surface $\Sigma$
\be
\label{k-6-rep-ter-ter}
\begin{array}{rcl}
\ds\mbox{div}^\Sigma\lf( \arctan\sigma\ \nabla^\Sigma\beta+\nabla^\Sigma\log\r\rg)&=&\ds{4}\,\frac{|(\nabla^\Sigma\r)^\perp|^2}{\r^2}+2\pi\, \sum_{p\in\r^{-1}(0)}\ \delta_{p}\ .
\end{array}
\ee
 Let $\ep>0$ and consider $\phi\in C^\infty_0(\Sigma\setminus \r^{-1}([0,\ep))$ and let $\delta>0$. 
Observe that $\arctan\sigma\ \nabla^\Sigma\beta+\nabla^\Sigma\log\r$ extends smoothly on $\Sigma\setminus \r^{-1}([0,\ep)$.
On $\Sigma$ we define 
\[
\chi_{\delta,\rho}:=\chi(\rho/\delta)
\]
where $\chi$ is the cut-off function on ${\R}_+$ defined by
\[
\chi(t)=\lf\{
\begin{array}{l}
1\quad\mbox{ for }t<1\\[3mm]
0\quad\mbox{ for }t>2
\end{array}
\rg.
\]
Observe that $\arctan\sigma\ \nabla^\Sigma\beta+\nabla^\Sigma\log\r$ extends continuously on $\Sigma\setminus \r^{-1}([0,\ep)$ and for any $\delta>0$ we write, using (\ref{k-6-rep}),
\be
\label{k-6-rep-quad}
\begin{array}{l}
\ds-\int_\Sigma \nabla^\Sigma\phi\ \lf( \arctan\sigma\ \nabla^\Sigma\beta+\nabla^\Sigma\log\r\rg)+{4}\,\phi\,\frac{|(\nabla^\Sigma\r)^\perp|^2}{\r^2}\ dvol_\Sigma\\[5mm]
\ds=-\int_\Sigma  \nabla^\Sigma(\chi_{\delta,\rho}\,\phi)\ \lf( \arctan\sigma\ \nabla^\Sigma\beta+\nabla^\Sigma\log\r\rg)+{4}\,\chi_{\delta,\rho}\,\phi\ \frac{|(\nabla^\Sigma\r)^\perp|^2}{\r^2}\ dvol_\Sigma
\end{array}
\ee
We have respectively, since $G$ is an immersion in ${\C}^2$
\be
\label{k-6-rep-quint}
\int_{\lf\{\rho<2\,\delta\ ;\ \r>\ep\rg\}\cap\,{Supp}(\phi)}\ dvol_\Sigma=O(\delta^2)\ \quad,\quad\ \|\nabla^\Sigma(\chi_{\delta,\rho}\,\phi)\|_\infty\le C_\phi\, \delta^{-1}
\ee
and
\be
\label{k-6-rep-sextuple}
\lf\| \arctan\sigma\ \nabla^\Sigma\beta+\nabla^\Sigma\log\r   \rg\|_{L^\infty(supp(\phi))}<C_\phi\quad,\quad \lf\|   \chi_{\delta,\rho}\,\phi\ \frac{|(\nabla^\Sigma\r)^\perp|^2}{\r^2} \rg\|_{L^\infty(supp(\phi))}<C_\phi\ ,
\ee
where $C_\phi$ is independent of $\delta$. Combining (\ref{k-6-rep-quad}), (\ref{k-6-rep-quint}) and (\ref{k-6-rep-sextuple}) and making $\delta$ converge to zero gives
\[
\int_\Sigma \nabla^\Sigma\phi\ \lf( \arctan\sigma\ \nabla^\Sigma\beta+\nabla^\Sigma\log\r\rg)+{4}\,\phi\,\frac{|(\nabla^\Sigma\r)^\perp|^2}{\r^2}\ dvol_\Sigma
\]
and this proves that  (\ref{k-6-rep-ter-ter}) holds on the whole surface $\Sigma$.

\medskip

Using again (\ref{kora-4re}) we obtain away from $\rho=0$
\be
\label{k-2}
\begin{array}{rcl}
\ds\r^3\, \nabla^\Sigma\beta\cdot\nabla^\Sigma\r&=&\rho^3\,\nabla^\Sigma\beta\cdot\nabla^\Sigma\rho+\nabla^\Sigma\varphi^2\cdot\nabla^\Sigma\beta\\[5mm]
&=&\ds 2^{-1}\,\rho^2\, \Delta^\Sigma\varphi+2\,\varphi\, \lf(1-4^{-1}\,\Delta^\Sigma\rho^2\rg)\\[5mm]
&=&\ds 2\,\varphi+2^{-1}\ div^\Sigma\lf(\rho^2\ \nabla^\Sigma\varphi-\varphi\,\nabla^\Sigma\rho^2\rg)\\[5mm]
&=&\ds 2\,\varphi+4^{-1}\ div^\Sigma\lf(\rho^4\ \nabla^\Sigma\sigma\rg)\\[5mm]
&=&\ds \sigma\, \rho^2+4^{-1}\ div^\Sigma\lf(\r^4\ \frac{\nabla^\Sigma\sigma}{1+\sigma^2}\rg)\\[5mm]
&=&\ds \r^2\, \frac{\sigma}{\sqrt{1+\sigma^2}}+4^{-1}\ div^\Sigma\lf(\r^4\ \frac{\nabla^\Sigma\sigma}{1+\sigma^2}\rg)\ .
\end{array}
\ee
Similarly as above for  (\ref{k-6-rep}) holds on the whole $\Sigma$, we prove that (\ref{k-2}) holds on $\Sigma$ away from $\r=0$.
Let $\chi$ be a smooth cut-off function on ${\R}_+$ defined above. Multiplying  (\ref{k-6-rep}) by $\chi(\r)$ and integrating over $\Sigma$ gives
\be
\label{k-7-rep}
\begin{array}{l}
\ds-\int_\Sigma\arctan\sigma\ \chi'(\r)\ \nabla^\Sigma\beta\cdot\nabla^\Sigma\r\ dvol_\Sigma-\int_\Sigma\, \chi'(\r)\ \frac{|\nabla^\Sigma\r|^2}{\r}\\[5mm]
\ds=4\,\int_\Sigma\chi(\r)\ \frac{|(\nabla^\Sigma\r)^\perp|^2}{\r^2}\ dvol_\Sigma+\theta_0
\end{array}
\ee
where $\theta_0=2\pi\, \mbox{Card}(\La^{-1}(\{0\})$.
Multiplying now (\ref{k-2}) by $-\,\arctan\sigma\ \chi'(\r)\, \r^{-3}$ and integrating over $\Sigma$ gives
\be
\label{k-8-rep}
\begin{array}{l}
\ds-\int_\Sigma\arctan\sigma\ \chi'(\r)\ \nabla^\Sigma\beta\cdot\nabla^\Sigma\r\ dvol_\Sigma=- \int_\Sigma\arctan\sigma\ \frac{\chi'(\r)}{\, \r}\ \frac{\sigma}{\sqrt{1+\sigma^2}}\ dvol_\Sigma
\\[5mm]
\ds\quad -\int_\Sigma \arctan\sigma\ \frac{\chi'(\r)}{4\, \r^3}\ div^\Sigma\lf(\r^4\ \frac{\nabla^\Sigma\sigma}{1+\sigma^2}\rg)\ dvol_\Sigma\ .
\end{array}
\ee
Integrating by parts gives for the second term of the r.-h.-s. of (\ref{k-8-rep})
\be
\label{k-9-rep}
\begin{array}{l}
\ds-\int_\Sigma \arctan\sigma\ \frac{\chi'(\r)}{4\, \r^3}\ div^\Sigma\lf(\r^4\ \frac{\nabla^\Sigma\sigma}{1+\sigma^2}\rg)\ dvol_\Sigma\\[5mm]
\ds=\int_\Sigma \frac{|\nabla^\Sigma\sigma|^2}{(1+\sigma^2)^2}\ \r\,\frac{\chi'(\r)}{4}\ \ dvol_\Sigma
\ds+\frac{1}{4}\,\int_\Sigma\arctan\sigma\ \chi''(\r)\ \r \,\nabla^\Sigma\r\cdot\frac{\nabla^\Sigma\sigma}{1+\sigma^2}\ dvol_\Sigma\\[5mm]
\ds\ -\frac{3}{4}\,\int_\Sigma\arctan\sigma\ \chi'(\r)\ \,\nabla^\Sigma\r\cdot\frac{\nabla^\Sigma\sigma}{1+\sigma^2}\ dvol_\Sigma
\end{array}
\ee
Combining (\ref{k-7-rep}), (\ref{k-8-rep}) and (\ref{k-9-rep}) gives then
\be
\label{k-10-rep}
\begin{array}{l}
\ds\int_\Sigma \frac{|\nabla^\Sigma\sigma|^2}{(1+\sigma^2)^2}\ \r\, \frac{\chi'(\r)}{4}\ \ dvol_\Sigma
\ds+\frac{1}{4}\,\int_\Sigma\arctan\sigma\ \chi''(\r)\ \r \,\nabla^\Sigma\r\cdot\frac{\nabla^\Sigma\sigma}{1+\sigma^2}\ dvol_\Sigma\\[5mm]
\ds\ -\frac{3}{4}\,\int_\Sigma\arctan\sigma\ \chi'(\r)\ \,\nabla^\Sigma\r\cdot\frac{\nabla^\Sigma\sigma}{1+\sigma^2}\ dvol_\Sigma\\[5mm]

\ds- \int_\Sigma\arctan\sigma\ \frac{\chi'(\r)}{\, \r}\ \frac{\sigma}{\sqrt{1+\sigma^2}}\ dvol_\Sigma
-\int_\Sigma\, \chi'(\r)\ \frac{|\nabla^\Sigma\r|^2}{\r}\ dvol_\Sigma\\[5mm]
\ds=4\,\int_\Sigma\chi(\r)\ \frac{|(\nabla^\Sigma\r)^\perp|^2}{\r^2}\ dvol_\Sigma+\theta_0
\end{array}
\ee
Observe that we have
\be
\label{k-11-rep}
\begin{array}{l}
\ds\frac{\nabla^\Sigma\sigma}{1+\sigma^2}=\frac{2}{\rho^2}\, \frac{\nabla^\Sigma\varphi}{1+\sigma^2}-4\,\frac{\varphi}{\rho^3}\ \frac{\nabla^\Sigma\rho}{1+\sigma^2}\\[5mm]
\ds =2\,\frac{\rho^2}{\r^4}\ \nabla^\Sigma\varphi-2\,\frac{\sigma}{\sqrt{1+\sigma^2}}\ \frac{\rho}{\r^2}\,\nabla^\Sigma\rho\ ,
\end{array}
\ee
and recall from (\ref{heat-13}) that
\[
|\nabla^\Sigma\rho|^2+\rho^{-2}\,|\nabla^\Sigma\varphi|^2=1\ .
\]
Hence we deduce that
\be
\label{k-12-rep}
\lf|\frac{\nabla^\Sigma\sigma}{1+\sigma^2}\rg|^2\le\,8\, \lf[\frac{\rho^6}{\r^8}+\frac{\rho^2}{\r^4}  \rg]\le \frac{16}{\r^2}\ .
\ee
Combining (\ref{k-10-rep}) and (\ref{k-12-rep}) gives the existence of a universal constant $C_0$ such that
\be
\label{k-13-rep}
2\pi\,\mbox{Card}\,{{\La}}^{\,-1}(0)+\int_{\r<1}\frac{|(\nabla^\Sigma\r)^\perp|^2}{\r^2}\ dvol_\Sigma\le C_0\ \int_{1<\r<2}\ dvol_\Sigma\ .
\ee
Recall
\be
\label{kora-36}
\r^3\,J(\nabla^H\r)=\rho^3\,J(\nabla^H\rho)+2\,\varphi\,J(\nabla^H\varphi)=-\rho^2\,\nabla^H\varphi+2\,\rho\, \varphi\, \nabla^H\rho=-\rho^4\,\nabla^H\lf(\frac{\varphi}{\rho^2}   \rg)
\ee
Since the immersion of $\Sigma$ by ${G}$ is lagrangian, the complex form $J$ is sending the tangential part of $\nabla^H\r$ (that we have denoted $\nabla^\Sigma\r$) to the normal directions to the surface and  vice versa that is  the normal part to the surface of $\nabla^H\r$ (that we have denoted $(\nabla^\Sigma\r)^\perp$) to the tangential directions.  Hence we have in particular
\be
\label{kora-36-a}
\r^3\, J(\nabla^\Sigma\r)^\perp=-\rho^4\,\nabla^\Sigma\lf(\frac{\varphi}{\rho^2}\rg)\ .
\ee
Hence
\be
\label{k-14-rep}
\frac{(\nabla^\Sigma\r)^\perp}{\r}=\frac{1}{2}\frac{\rho^4}{\r^4}\ J(\nabla^\Sigma\sigma)=\frac{1}{2}\, J\lf(  \frac{\nabla^\Sigma\sigma}{1+\sigma^2}  \rg)\ .
\ee
Combining (\ref{k-13-rep}) and (\ref{k-14-rep}) gives the  lemma~\ref{lm-mono}.\hfill$\Box$

\section{Proof of the main theorem~\ref{th-main}.}
Let $0<r<1$. Observe that the lower bound in (\ref{000}) is a direct consequence of (\ref{k-15-rep}) after rescaling at $r$. We now prove the upper bound in (\ref{000}). 

\medskip

We  consider a smooth cut-off function $\chi$ on ${\R}_+$ such that
\[
\chi(t)=\lf\{
\begin{array}{l}
1\quad\mbox{ for }t<1\\[3mm]
0\quad\mbox{ for }t>2
\end{array}
\rg.\quad,\quad \chi'\le 0\quad\mbox{on }{\R}_+\quad\mbox{ and }\quad\chi'\equiv-3/2\quad\mbox{ on }\quad[5/4,7/4]\quad.
\]
Replacing in (\ref{k-10-rep}) $\chi$ by $\chi(\r/r)$ and shifting the first term of the l.h.s to the r.h.s gives
\be
\label{k-16-rep}
\begin{array}{l}
\ds\frac{1}{4}\,\int_\Sigma\arctan\sigma\ r^{-2}\ \chi''(\r/r)\ \r \,\nabla^\Sigma\r\cdot\frac{\nabla^\Sigma\sigma}{1+\sigma^2}\ dvol_\Sigma\\[5mm]
\ds-\frac{3}{4}\,\int_\Sigma\arctan\sigma\ r^{-1}\ \chi'(\r/r)\ \,\nabla^\Sigma\r\cdot\frac{\nabla^\Sigma\sigma}{1+\sigma^2}\ dvol_\Sigma\\[5mm]
\ds-\int_\Sigma \frac{\chi'(\r/r)}{r\, \r}\ \lf[\frac{\sigma\,\arctan\sigma}{\sqrt{1+\sigma^2}}+ {|\nabla^\Sigma\r|^2}\rg]\ dvol_\Sigma\\[5mm]
\ds=4\,\int_\Sigma\chi(\r/r)\ \frac{|(\nabla^\Sigma\r)^\perp|^2}{\r^2}\ dvol_\Sigma-\, \int_\Sigma \frac{|\nabla^\Sigma\sigma|^2}{(1+\sigma^2)^2}\ \frac{\r}{r}\,\frac{\chi'(\r/r)}{4}\ \ dvol_\Sigma+\theta_0
\end{array}
\ee
Recall from  (\ref{k-14-rep}) that 
\[
\frac{1}{\sqrt{1+\sigma^2}}=|\nabla^H\r|^2=|\nabla^\Sigma\r|^2+|(\nabla^\Sigma\r)^\perp|^2\ ,
\]
Hence (\ref{k-16-rep}) becomes
\be
\label{k-17-rep}
\begin{array}{l}
\ds\frac{1}{4}\,\int_\Sigma\arctan\sigma\ r^{-2}\ \chi''(\r/r)\ \r \,\nabla^\Sigma\r\cdot\frac{\nabla^\Sigma\sigma}{1+\sigma^2}\ dvol_\Sigma\\[5mm]
\ds-\frac{3}{4}\,\int_\Sigma\arctan\sigma\ r^{-1}\ \chi'(\r/r)\ \,\nabla^\Sigma\r\cdot\frac{\nabla^\Sigma\sigma}{1+\sigma^2}\ dvol_\Sigma\\[5mm]
\ds-\int_\Sigma \frac{\chi'(\r/r)}{r\, \r}\ \frac{\sigma\,\arctan\sigma+1}{\sqrt{1+\sigma^2}}\ dvol_\Sigma\\[5mm]
\ds=4\,\int_\Sigma\chi(\r/r)\ \frac{|(\nabla^\Sigma\r)^\perp|^2}{\r^2}\ dvol_\Sigma-\int_\Sigma {\chi'(\r/r)}\ \frac{|(\nabla^\Sigma\r)^\perp|^2}{r\,\r}\ dvol_\Sigma\\[5mm]
\ds\quad-\, \int_\Sigma \frac{|\nabla^\Sigma\sigma|^2}{(1+\sigma^2)^2}\ \frac{\r}{r}\,\frac{\chi'(\r/r)}{4}\ \ dvol_\Sigma+\theta_0
\end{array}
\ee
Observe that for  any $\sigma\in {\R}$
\be
\label{k-18-rep}
\lf(\frac{\sigma\,\arctan\sigma+1}{\sqrt{1+\sigma^2}}\rg)'= \frac{\arctan\sigma}{(1+\sigma^2)^{3/2}}\ .
\ee
We deduce that
\be
\label{k-19-rep}
\forall \,\sigma\in {\R}\quad\quad\quad1\le \frac{\sigma\,\arctan\sigma+1}{\sqrt{1+\sigma^2}}\le \frac{\pi}{2}\ .
\ee
This gives in particular the existence of a universal constant $C_1$ for $r<1/2$, using the fact that $\chi'\le 0$ and $\chi'=-3/2$ on $[5/4,7/4]$
\be
\label{k-20-rep}
\begin{array}{l}
\ds r^{-2}\,\int_{ 5\, r/4<\r<7\, r/4}\ dvol_\Sigma\le\ C_1\ \int_{\r<2\,r}\lf|\frac{\nabla^\Sigma\sigma}{1+\sigma^2}\rg|^2\ dvol_\Sigma+C_1\ \theta_0\\[5mm]
\ds+ C_1\, r^{-1}\, \int_{r<\r<2 r}\lf|\frac{\nabla^\Sigma\sigma}{1+\sigma^2}\rg|\ dvol_\Sigma
\end{array}
\ee
Using Cauchy Schwartz together with (\ref{k-15-rep}) we obtain
\be
\label{k-21-rep}
\begin{array}{l}
\ds r^{-2}\,\int_{ 5\, r/4<\r<7\, r/4}\ dvol_\Sigma\le\ C_1\ \int_{ 1<\r<2}\ dvol_\Sigma\\[5mm]
\ds\quad+ C_2\ \lf[ \int_{ 1<\r<2}\ dvol_\Sigma \rg]^{1/2}\ \lf[r^{-2}\ \int_{ r<\r<2\,r}\ dvol_\Sigma\rg]^{1/2}
\end{array}
\ee
where the $C_i$s denote universal constants. Let 
\[
A:=\sup_{r<1/2}r^{-2}\,\int_{ 5\, r/4<\r<7\, r/4}\ dvol_\Sigma
\]
We deduce from (\ref{k-21-rep})
\be
\label{k-22-rep}
A\le C_3 \int_{ 1/2<\r<2}\ dvol_\Sigma+ C_3\  \lf[ \int_{ 1/2<\r<2}\ dvol_\Sigma \rg]^{1/2} A^{1/2}
\ee
This last inequality implies the upper bound in (\ref{000}) after observing that 
\be
\label{k-23-rep}
 \frac{1}{r^2}\int_{\r<r}\ dvol_\Sigma=\frac{1}{r^2}\lf[ \sum_{j=0}^\infty\int_{4^{-j-1}\,r<\r<4^{-j}\,r} \ dvol_\Sigma \rg]\le\lf[ \sum_{j=0}^\infty 2\, 4^{-2j}\, A \rg]\le C A
\ee
This concludes the proof of theorem~\ref{th-main}.\hfill $\Box$.

\medskip

Observe that a similar approach holds for the proof of the classical monotonicity formula (\ref{001}). Indeed, the minimal surface equation in ${\R}^n$ is $\Delta^\Sigma\Phi=0$ hence we deduce
\be
\label{k-24-rep}
\Delta^\Sigma\rho^2=2\,|\nabla^\Sigma\Phi|^2=4
\ee
away from  $\rho=0$ this last equation implies
\be
\label{k-25-rep}
\Delta^\Sigma\log\rho^2=\frac{4}{\rho^2}-4\,\frac{|\nabla^\Sigma\rho|^2}{\rho^2}
\ee
Since with the above notations $1=|\nabla\rho|^2=|\nabla^\Sigma\rho|^2+|(\nabla^\Sigma\rho)^\perp|^2$, the previous identity implies as when passing from (\ref{k-6-rep}) to (\ref{k-6-rep-ter-ter})
\be
\label{k-26-rep}
\frac{1}{4}\,\Delta^\Sigma\log\rho^2=\frac{|(\nabla^\Sigma\rho)^\perp|^2}{\rho^2}+\pi\,\sum_{p\in\rho^{-1}(\{0\})}\delta_p
\ee
Integrating over $\rho<r$ and using again (\ref{k-24-rep}) gives
\be
\label{k-27-rep}
\begin{array}{l}
\ds\frac{1}{r^2}\int_{\rho<r}\ dvol_\Sigma= \frac{1}{4\,r^2}\int_{\rho<r}\Delta^\Sigma\rho^2\ dvol_\Sigma=\frac{1}{4}\,\int_{\rho=r}\ \frac{\p_\nu\rho^2}{\rho^2}\ dl_\Sigma\\[5mm]
\ds=\int_{\rho<r}\frac{|(\nabla\rho)^\perp|^2}{\rho^2}\ dvol_\Sigma+\pi\,  \mbox{Card}\lf( \Phi^{-1}(0)\rg)
\end{array}
\ee
This implies the monotonicity formula (\ref{001}).

\section{Proof of the Bernstein-Liouville theorem~\ref{th-berns-liouv}}
Modulo a translation (by the Heisenberg group action) we can assume that the surface is passing though the origin and that, using the above notations, we have
\be
\label{IV.1}
\theta_0\ge 2\pi\ .
\ee
Since by assumption
\be
\label{IV.1bis}
\limsup_{\r\rightarrow +\infty}\frac{1}{r^2}\int_{r<\r<2\,r}\ dvol_{\Sigma}<+\infty
\ee
we deduce from  lemma~\ref{lm-mono}
\[
\int_\Sigma\lf|\frac{\nabla^\Sigma\sigma}{1+\sigma^2}\rg|^2\ dvol_\Sigma<+\infty
\]
This implies in particular that
\be
\label{IV.2}
\lim_{r\rightarrow +\infty}\int_{r<\r<2\,r}\lf|\frac{\nabla^\Sigma\sigma}{1+\sigma^2}\rg|^2\ dvol_{\Sigma}=0\ .
\ee 
We fix $\ep>0$ and we  consider a smooth cut-off function $\chi$ on ${\R}_+$ such that
\[
\chi_\ep(t)=\lf\{
\begin{array}{l}
1\quad\mbox{ for }t<1\\[3mm]
0\quad\mbox{ for }t>2
\end{array}
\rg.\quad,\quad \chi_\ep'\le 0\quad\mbox{on }{\R}_+\quad\mbox{ and }\quad\chi_\ep'\equiv-1\quad\mbox{ on }\quad[1+\ep,2-\ep]\quad.
\]
and we can also assume that $\|\chi_\ep'\|_\infty$  is uniformly bounded independently of $\ep$. We shall denote $C_\ep:=\|\chi_\ep''\|_\infty$ (which clearly is not bounded independently of $\ep$). Replacing in (\ref{k-10-rep}) $\chi$ by $\chi_\ep(\r/r)$, using (\ref{k-14-rep}) and shifting the first term of the l.h.s to the r.h.s gives
\be
\label{IV.3}
\begin{array}{l}
\ds\frac{1}{4}\,\int_\Sigma\arctan\sigma\ r^{-2}\ \chi_\ep''(\r/r)\ \r \,\nabla^\Sigma\r\cdot\frac{\nabla^\Sigma\sigma}{1+\sigma^2}\ dvol_\Sigma\\[5mm]
\ds-\frac{3}{4}\,\int_\Sigma\arctan\sigma\ r^{-1}\ \chi_\ep'(\r/r)\ \,\nabla^\Sigma\r\cdot\frac{\nabla^\Sigma\sigma}{1+\sigma^2}\ dvol_\Sigma\\[5mm]
\ds-\int_\Sigma \frac{\chi_\ep'(\r/r)}{r\, \r}\ \frac{\sigma\,\arctan\sigma+1}{\sqrt{1+\sigma^2}}\ dvol_\Sigma\\[5mm]
\ds=\,\int_\Sigma\chi_\ep(\r/r)\ \frac{|\nabla^\Sigma\sigma|^2}{(1+\sigma^2)^2}\ dvol_\Sigma-\int_\Sigma {\chi_\ep'(\r/r)}\ \frac{\r}{4\,r}\ \frac{|\nabla^\Sigma\sigma|^2}{(1+\sigma^2)^2}dvol_\Sigma\\[5mm]
\ds\quad-\, \int_\Sigma \frac{|\nabla^\Sigma\sigma|^2}{(1+\sigma^2)^2}\ \frac{\r}{r}\,\frac{\chi_\ep'(\r/r)}{4}\ \ dvol_\Sigma+\theta_0
\end{array}
\ee
We have respectively for $\ep>0$ fixed
\be
\label{IV.4}
\begin{array}{l}
\ds\lf|\frac{1}{4}\,\int_\Sigma\arctan\sigma\ r^{-2}\ \chi_\ep''(\r/r)\ \r \,\nabla^\Sigma\r\cdot\frac{\nabla^\Sigma\sigma}{1+\sigma^2}\ dvol_\Sigma\rg|\\[5mm]
\ds\ \le C_\ep\lf[  \int_{r<\r<2\,r}\frac{\r^2}{r^4}\ dvol_\Sigma\rg]^{1/2}\lf[  \int_{r<\r<2\,r}\lf|\frac{\nabla^\Sigma\sigma}{1+\sigma^2}\rg|^2\ dvol_{\Sigma} \rg]^{1/2}\rightarrow 0\quad\mbox{ as }r\rightarrow+\infty
\end{array}
\ee
and
\be
\label{IV.4-b}
\begin{array}{l}
\ds\lf|  \frac{3}{4}\,\int_\Sigma\arctan\sigma\ r^{-1}\ \chi_\ep'(\r/r)\ \,\nabla^\Sigma\r\cdot\frac{\nabla^\Sigma\sigma}{1+\sigma^2}\ dvol_\Sigma\rg|\\[5mm]
\ds\ \le\|\chi'_\ep\|_\infty\lf[  \int_{r<\r<2\,r}\frac{1}{r^2}\ dvol_\Sigma\rg]^{1/2}\lf[  \int_{r<\r<2\,r}\lf|\frac{\nabla^\Sigma\sigma}{1+\sigma^2}\rg|^2\ dvol_{\Sigma} \rg]^{1/2}\rightarrow 0\quad\mbox{ as }r\rightarrow+\infty
\end{array}
\ee
From the assumption of the theorem we have that $\sigma$ converges uniformly towards 0 as $\r\rightarrow+\infty$ hence
\be
\label{IV.5}
\lim_{r\rightarrow +\infty}-\int_\Sigma \frac{\chi_\ep'(\r/r)}{r\, \r}\ \frac{\sigma\,\arctan\sigma+1}{\sqrt{1+\sigma^2}}\ dvol_\Sigma+\int_{r<\r<2\,r} \frac{\chi_\ep'(\r/r)}{r\, \r}\ \ dvol_\Sigma=0
\ee
Because of (\ref{IV.1bis}) we have the existence of $r_k\rightarrow +\infty$ such that
\be
\label{IV.6}
\frac{1}{r_k^2}\int_{r_k<\r<r_k+\ep}\ \ dvol_\Sigma+\frac{1}{r_k^2}\int_{2\,r_k-\ep<\r<2\,r_k}\ \ dvol_\Sigma\le C\ \ep
\ee
Hence, since $\|\chi'_\ep\|_\infty$ is bounded independently of $\ep$, we deduce from (\ref{IV.5}), (\ref{IV.6})
\be
\label{IV.7}
\lim_{\ep\rightarrow 0}\limsup_{k\rightarrow +\infty}\lf|-\int_\Sigma \frac{\chi_\ep'(\r/r_k)}{r_k\, \r}\ \frac{\sigma\,\arctan\sigma+1}{\sqrt{1+\sigma^2}}\ dvol_\Sigma-\int_{r_k<\r<2\,r_k} \frac{1}{r_k\, \r}\ \ dvol_\Sigma\rg|=0
\ee
Which implies by assumption
\be
\label{IV.8}
\lim_{\ep\rightarrow 0}\limsup_{k\rightarrow +\infty}\lf|-\int_\Sigma \frac{\chi_\ep'(\r/r_k)}{r_k\, \r}\ \frac{\sigma\,\arctan\sigma+1}{\sqrt{1+\sigma^2}}\ dvol_\Sigma-2\pi\rg|=0
\ee
Combining now the fact that $\theta_0\ge 2\pi$ together with (\ref{IV.3}), (\ref{IV.4}), (\ref{IV.4-b}) and (\ref{IV.8}) we obtain
\be
\label{IV.9}
0=\int_\Sigma\lf|\frac{\nabla^\Sigma\sigma}{1+\sigma^2}\rg|^2\ dvol_\Sigma=\int_\Sigma \frac{|(\nabla^\Sigma\r)^\perp|^2}{\r^2}\ dvol_\Sigma\ .
\ee
Hence  $\nabla^\Sigma\r=\nabla^H\r$ and $\La$ is a H-minimal legendrian smooth conical immersion. The only smooth conical H-minimal legendrian immersions are lagrangian planes (see \cite{SW2}). This concludes
the proof of theorem~\ref{th-berns-liouv}.\hfill $\Box$


\begin{thebibliography}{99}
\bibitem{Oh1}  Oh, Yong-Geun Second variation and stabilities of minimal Lagrangian submanifolds in K\"ahler manifolds. Invent. Math. 101 (1990), no. 2, 501-519.
\bibitem{Oh2}  Oh, Yong-Geun Volume minimization of lagrangian submanifolds under hamiltonian deformations, Math Z. 212 (1993) 175-192.
 \bibitem{Piga-Riv-Lag-Visc} Pigati, Alessandro; Rivi\`ere, Tristan ``The Viscosity Method for Minimal Lagrangian Surfaces in K\"ahler-Einstein Surfaces'' in preparation.
 \bibitem{Riv} Rivi\`ere, Tristan   ``Minmax Hierarchies, Minimal Fibrations and a PDE based Proof of the Willmore Conjecture''    https://arxiv.org/abs/2007.05467
 \bibitem{SW1} Schoen, Richard; Wolfson, Jon Minimizing volume among Lagrangian submanifolds. Differential equations: La Pietra 1996 (Florence), 181-199, Proc. Sympos. Pure Math., 65, Amer. Math. Soc., Providence, RI, 1999.
 \bibitem{SW2} Schoen, R.; Wolfson, J. Minimizing area among Lagrangian surfaces: the mapping problem. J. Differential Geom. 58 (2001), no. 1, 1-86.
\bibitem{Sim} Simon Leon, Lectures on Geometric Measure Theory, Australian National University Centre for Mathematical Analysis, Canberra, 1983
 \end{thebibliography}
\end{document}